\documentclass[10pt]{amsart}
\begin{document}
\title{Closures of $O_n$-orbits in the flag variety for $GL_n$, II}
\author{William M. McGovern}
\begin{abstract}
We give a necessary and sufficient pattern avoidance criterion for the conjugates of the bottom vertex of the Bruhat graph attached to an $O_n$ orbit $\mathcal O$ in the flag variety of $GL_n$ to have degree equal to the rank of this graph as a poset, showing also that this condition is equivalent to the rational smoothness of $\overline{\mathcal O}$.  We also give a necessary and sufficient condition in terms of pattern avoidance for $\overline{\mathcal O}$ to be smooth.
\end{abstract}
\maketitle

This short note is a continuation of \cite{M19}, whose notation and references we shall preserve; there is an additional reference, namely \cite{MS05}.  Our main result shows that the graph-theoretic and pattern avoidance criteria for rational smoothness of $O_n$ orbit closures in the flag variety for $GL_n$ is both necessary and sufficient in all cases.

\newtheorem*{theorem1}{Theorem 1}
\begin{theorem1}
For all $n$, the orbit closure $\overline{\mathcal O}_\pi$ is rationally smooth if and only if all conjugates of $w_0$ in the order ideal $I_\pi$ have degree $r(\pi)$, or if and only if $\pi$ avoids all bad patterns of \cite{M19}.
\end{theorem1}

\begin{proof}
We begin with a small clarification of the proof of Theorem 2 in \cite{M19}:  the slice $\mathcal S$ of the orbit $\overline{\mathcal O}$  to $\mathcal O_c$ at $P$ is defined by the vanishing of all minors in the upper left $i\times j$ submatrix of the Gram matrix $G_\pi$ of size $(\pi_{ij}+1\times \pi_{ij}+1)$, as $i,j$ runs through all indices at most equal to $n$ \cite{MS05}; these minors generate a prime ideal.   It was shown in \cite{M19} that if one of the bad patterns is included, then some conjugate of $w_0$ in $I_\pi$ has degree larger than $r(\pi)$.  If all bad patterns are avoided, then repeating the argument in the proof of Theorem 2 and using in particular that the patterns 2137654 and 21435 (the latter in effect a bad pattern for all involutions of odd length) are avoided, we find that the vanishing minors used to define the slice $\mathcal S''$ are all distinct and that this slice can be further sliced to obtain a variety rationally smooth away from the origin, so that all three parts of Brion's criterion hold and $\overline{\mathcal O}_\pi$ is rationally smooth at the point $P$ defined in the proof of Theorem 2 and so everywhere.
\end{proof}

As shown in Theorem 3 of \cite{M19}, we can improve this result for $n$ even; in that case $\overline{\mathcal O}_\pi$ is rationally smooth if and only if the degree of just $w_0$ in $I_\pi$ is $r(\pi)$. 

Our other result gives a necessary and sufficient condition for smoothness of orbit closures in terms of pattern avoidance; there is no graph-theoretic analogue of this condition.  We find that two additional bad patterns must be avoided, namely 2143 (in all cases, regardless of the number of fixed points between the 21 and the 43) and 1324:

\newtheorem*{theorem2}{Theorem 2}
\begin{theorem2}
For any involution $\pi$ containing one of the twenty-four bad patterns of \cite{M19} or the pattern 2143 or 1324, the corresponding orbit closure $\overline{\mathcal O}_\pi$ is not smooth.  Conversely, if $\pi$ avoids all of these patterns then $\overline{\mathcal O}_\pi$ is smooth.
\end{theorem2}

\begin{proof}
This follows from \cite{M19} if one of the twenty-four bad patterns is included.  If all are avoided but 2143 or 1324 is included, so that $\overline{\mathcal O}_\pi$ is rationally smooth, then as in the proof of Theorem 2 we define the slice $\mathcal S$ by equating to 0 all $(\pi_{ij}+1)\times(\pi_{ij}+1)$ minors in all upper left $i\times j$ submatrices of the matrix $G$ in the proof of Theorem 2, as $i,j$ run through all indices between $1$ and $n$.  These minors generate the ideal $J_\pi$ of definition of the orbit closure $\overline{\mathcal O}_\pi$ \cite{MS05}.  Arguing inductively, as the proof of Theorem 1, starting with each of the bad patterns and adding fixed points and pairs of flipped indices one at a time, we see that by equating these minors to 0 (as we did to define the slice $\mathcal S$ in the proof of Theorem 2) we can solve for some variables in terms of the others; all terms of the equations in the remaining variables have degree at least 2.  Thus smoothness fails at the point $P$ defined in the proof of Theorem 2, where all variables are 0.  Conversely, if all bad patterns are avoided, then the generators of $J_\pi$ all include terms of degree 1, whence $\overline{\mathcal O}_\pi$ is smooth, as desired.
\end{proof}

\end{document}